\documentclass[reqno,12pt]{amsart}
\usepackage{latexsym, amsfonts,mathrsfs, amsmath, amssymb, amscd, epsfig}
\usepackage{mathrsfs}
\usepackage{cases}
\usepackage{color}

\newtheorem{thm}{Theorem}

\newtheorem{prop}[thm]{Proposition}
\newtheorem{lem}[thm]{Lemma}

\newtheorem{rmk}{Remark}

\newenvironment{pf}{{\noindent \it \bf Proof:}}{{\hfill$\Box$}\\}
\newenvironment{pf2}{{\noindent \it \bf Proof of Theorem \ref{thm2}:}}{{\hfill$\Box$}\\}
\def\Div{{\rm div}}

\def\mcC{\mathcal{C}}

\def\mcA{\mathcal{A}}

\def\bu{{\mathbf u}}
\def\bv{{\mathbf v}}
\def\bw{{\mathbf w}}
\def\bn{{\mathbf n}}

\begin{document}

\title[ ]{Three Dimensional Steady Subsonic Euler Flows in Bounded Nozzles}

\author[ ]{Chao Chen}
\address{School of Mathematics and Computer Science, Fujian Normal University, Fujian, 350108, China}
\email{chenchao\_math@sina.cn}
\author[ ]{Chunjing Xie}
\address{Department of mathematics, Institute of Natural Sciences, Ministry of Education Key Laboratory of Scientific and Engineering Computing, Shanghai Jiao Tong University, 800 Dongchuan Road, Shanghai, 200240, China}
\email{cjxie@sjtu.edu.cn}

%\date{ }

\begin{abstract} In this paper, we study the existence and uniqueness of three dimensional steady Euler flows in rectangular nozzles when prescribing normal component of momentum at both the entrance and exit. If, in addition, the normal component of the voriticity  and the variation of Bernoulli's function at the exit are both zero, then there exists a unique subsonic potential flow when the magnitude of the normal component of the momentum is less than a critical number. As the magnitude of the normal component of the momentum approaches the critical number, the associated flows converge to a subsonic-sonic flow. Furthermore, when the normal component of vorticity and the variation of Bernoulli's function are both small, the existence of subsonic Euler flows is established. The proof of these results is based on a new formulation for the Euler system,  a priori estimate for nonlinear elliptic equations with nonlinear boundary conditions, detailed study for  a linear div-curl system, and delicate estimate for the transport equations.
\end{abstract}

\maketitle

\section{Introduction and Main Results}\label{Formulation}

Multidimensional gas flows give rise many outstanding challenging
problems. The steady fluid is a natural starting point for the study on the multidimensional flows. However, the steady Euler equations themselves are
not easy to tackle, since the equations may not only be hyperbolic
or hyperbolic-elliptic coupled system, but also have discontinuous
solutions such as shock waves and vortex sheets.  An important
approximate model is the potential flow, which describes flows without vorticity. Since 1950's, tremendous progress
has been made on the study for potential flows. Subsonic potential
flows around a body were studied extensively by
Shiffman\cite{Shiffman}, Bers\cite{Bers1,Bers2}, Finn,
Gilbarg\cite{FG1,FG2}, and Dong\cite{Dong}, et al.  The existence of subsonic potential flows in  multidimensional infinitely long nozzles was achieved in
 \cite{XX1, XX2, DXY}. The subsonic-sonic
flow as a limit of subsonic flows were studied in
\cite{CDSW, XX1, XX2,Huang} via compensated compactness method. Subsonic flows with a sonic boundary were constructed in \cite{WX}.

When the flows have non-zero vorticity, steady incompressible Euler flows which are the zero Mach number limits of compressible subsonic flows were investigated in \cite{Alber, Glass, TX} and references therein.
The existence of subsonic Euler flows in two dimensional and three dimensional axially symmetric nozzles was established in \cite{XX3, DD} via a stream function formulation for the Euler equations, see \cite{ChenJun, ChenDX} for recent study on subsonic Euler flows in half plane.

As an important constituent of transonic flows, the existence of subsonic flows which are small perturbations of certain simple background flows was also studied in
\cite{LXY1,LXY2,ChenS, ChenY,XYY,XY3, Yuan} and references
therein.

This paper studies general subsonic Euler flows when prescribing the normal component of the momentum on the boundary.
The three dimensional steady isentropic compressible ideal flows are governed by
the following Euler equations:
\begin{eqnarray}
\nabla\cdot(\rho\mathbf{u})=0,\label{Massconservation}\\
\nabla\cdot(\rho\mathbf{u}\otimes\mathbf{u})+\nabla
p=0,\label{Momentumconservation}
\end{eqnarray}
where $\rho$, $\mathbf{u}=(u_1,u_2,u_3)$, $p$ are the
density, velocity, and pressure, respectively. We assume that $p=p(\rho)$ satisfies
\begin{equation}
p'(\rho)>0 \quad \text{and}\quad p''(\rho)>0\quad \text{for} \quad \rho>0.
\end{equation}
 The equations (\ref{Massconservation}) and (\ref{Momentumconservation}) form a hyperbolic-elliptic coupled system for subsonic flows
($|\mathbf{u}|^2<c^2$) and a hyperbolic system for supersonic
flow ($|\mathbf{u}|^2>c^2$) (cf. \cite{XX3}), where the quantity  $c(\rho)=\sqrt{p'(\rho)}$ is called local sound speed.

In this paper, we study steady Euler flows in a rectangular domain $\Omega=[0, L]\times[0,1]^2$ with
the entrance $\Gamma_{-}= \{0\}\times [0,1]^2$, the exit $\Gamma_{+}=\{L\}\times [0,1]^2$, and the solid wall $\Gamma =\partial \Omega\backslash (\Gamma_-\cup \Gamma_+)$.

 We prescribe the normal component of the momentum on the boundary, i.e.,
\begin{equation}\label{Boundary}
\rho\mathbf{u}\cdot\mathbf{n}=f,
\end{equation}
where $\mathbf{n}$ is the unit outer normal vector and $f$ satisfies the compatibility
condition
\begin{equation}\label{compatibility1}
\int_{\partial\Omega}f dS=0.
\end{equation}
Furthermore, at the entrance $\Gamma_-$ and the exit $\Gamma_+$, $f$ satisfies
\begin{equation}\label{positive}
f<0\,\,\ \text{on}\,\ \Gamma_-;\quad f=0\,\,\ \text{on}\,\ \Gamma
;\quad f>0\,\,\ \text{on}\,\ \Gamma_+ .
\end{equation}

Since three dimensional steady subsonic Euler system has two elliptic modes and two hyperbolic modes, one needs also to impose conditions for hyperbolic modes at the entrance. It follows from the Euler equations that
\begin{equation}
\bu\cdot \nabla B=0
\end{equation}
where $B=\frac{1}{2}|\bu|^2+h(\rho)$ is called Bernoulli function with the enthalpy $h$ defined by $h(\rho)=\int_a^\rho\frac{p'(s)}{s}d s$ for some constant $a$.
At the entrance $\Gamma_-$, the following boundary conditions are prescribed
\begin{equation}
(\nabla\times\bu)\cdot\mathbf{n} = \kappa,\,\,\,\, B=B_0\,\
\text{on}\,\ \Gamma_-.\label{Entrance2}
\end{equation}

Before we state our main results, let us introduce some notations in this paper. $\|\cdot\|_{L^p(\Omega)}$ and $\|\cdot\|_{H^k(\Omega)}$ are standard norms for $L^p$ spaces and Sobolev spaces $H^k(\Omega)$. Let  $C^{k,\alpha}(\Omega)$ ($C^{k, \alpha}(\bar\Omega)$) be the standard H\"{o}lder spaces and denote
$\|u\|_{2, \alpha; \Omega}=\|u\|_{C^{2, \alpha}(\bar\Omega)}$.
Define $\mcA = \{u|u\in L^1(\Omega), \int_\Omega udx=0\}$ and  $\mcC^{k, \alpha}(\Omega)=C^{k, \alpha}(\Omega)\cap \mcA$ ($\mcC^{k, \alpha}(\bar\Omega)=C^{k, \alpha}(\bar\Omega)\cap \mcA$). Let $\Gamma_2=\{x=(x_1,x_2,x_3)\in\partial\Omega|x_2=0 \text{ or } 1\}$ and $\Gamma_3=\{x=(x_1,x_2,x_3)\in\partial\Omega|x_3=0\text{ or } 1\}$ and denote $\Gamma=\Gamma_2\cup\Gamma_3$. Then $\mathcal{E}_e=( \Gamma\cap\Gamma_-)\cup(\Gamma\cap\Gamma_+)$ and $\mathcal{E}_{\Gamma}=\Gamma_2\cap\Gamma_3$ is the edge of the nozzle boundary and $\mathscr{C}=\mathcal{E}_e\cap\mathcal{E}_{\Gamma}$ is the set of the corner points on the nozzle boundary.

Our first result is the following:
\begin{thm}\label{thm1}
If $\kappa\equiv 0$ and $B_0\equiv \bar B$(constant), then for any given $f$ satisfying \eqref{compatibility1} and
\begin{equation}\label{comp2}
\frac{\partial f}{\partial \nu}=0\quad \text{on}\quad \partial\Gamma_-\cup\partial\Gamma_+,
\end{equation}
where $\nu$ is the unit normal of $\partial\Gamma_-$ ($\partial\Gamma_+$) on the plane $\Gamma_-$ ($\Gamma_+$), there exists $\theta^*>0$ such that for any $\theta\in (0, \theta^*)$, the problem \eqref{Massconservation} and \eqref{Momentumconservation} supplemented with boundary conditions \eqref{Entrance2}
and
\begin{equation}\label{Boundary2}
\rho\mathbf{u}\cdot\mathbf{n}=\theta f,
\end{equation}
has a unique subsonic solution $(\rho, \bu)$
which satisfies $\inf_{\bar\Omega}u_1>0$. Furthermore, as $\theta\to \theta^*$, the maximum of flow Mach numbers, i.e., $\max_{\bar\Omega}\frac{|\bu|}{c(\rho)}$, goes to $1$. Finally, as $\theta\to \theta^*$, the associated solutions $(\rho, \bu)$ have an almost everywhere convergent subsequence, whose limit satisfies the Euler system
\begin{equation}
\left\{
\begin{aligned}
&\Div \left(h^{-1}\left(\bar B-\frac{|\bu|^2}{2}\right)\bu\right)=0,\\
&\nabla\times\bu =0
\end{aligned}
\right.
\end{equation}
in the sense of distribution and the boundary condition \eqref{Boundary2} with $\theta=\theta^*$ in the sense of boundary trace.
\end{thm}

For the general case where $\kappa$ and $B_0$ in \eqref{Entrance2} are not zero and constant, respectively, we have the following result.
\begin{thm}\label{thm2}
 If the problem \eqref{Massconservation} and \eqref{Momentumconservation} with boundary conditions \eqref{Boundary} and \eqref{Entrance2} with $f=\bar{f}, \kappa\equiv 0, B_0\equiv\bar{B}$ has a unique subsonic solution $(\bar{\rho},\bar \bu)$ satisfying $\inf_{\bar\Omega}\bar{u}_1>0$ in $\Omega$, then there exists an $\epsilon_0>0$ such that if $\kappa\in C^{1,\alpha}(\Gamma_-)$, $B_0\in C^{2, \alpha}(\Gamma_-)$ satisfy
\begin{equation}\label{B0}
B_0=\bar{B},\quad\kappa=0  \quad\text{and}\,\,\,\,\frac{\partial B_0}{\partial \nu}=0\quad \text{on}\quad \partial\Gamma_-
\end{equation}
 where $\nu$ is  the unit normal of $\partial\Gamma_-$ on $\Gamma_-$, and
\begin{equation}\label{*}
\|f-\bar{f}\|_{C^{2,\alpha}(\bar{\Omega})}+\|\kappa\|_{C^{1,\alpha}(\Gamma_{-})}+\|B_0-\bar{B}\|_{C^{2,\alpha}(\Gamma_-)}\leq\epsilon_0,
\end{equation}
then there exists  a $C^{2,\alpha}$-smooth
solution to the original Euler system
(\ref{Massconservation})-(\ref{Momentumconservation}) with the
boundary conditions (\ref{Boundary}) and (\ref{Entrance2}).
\end{thm}

\begin{rmk}
The results in this paper can be generalized to domains with more complicated geometry, see our forthcoming paper \cite{ChenXie2}.
\end{rmk}

The paper is organized as follows. In section \ref{secpreliminary} we give an equivalent formulation for 3D Euler system and study a linear div-curl system which used in Section \ref{Existence}. Section \ref{PotentialFlows} devotes to study the potential flows and gives the proof for Theorem \ref{thm1}. Theorem 2 is proved in \ref{Existence}. There are two appendices at the end of the paper, where the estimate along the streamlines and the regularity estimate for the elliptic equations with compatibility conditions are investigated.

%\newpage

\section{Preliminaries}\label{secpreliminary}

In this section, we first derive an equivalent formulation for the Euler system. Second, a linear div-curl
system appeared in the proof of Theorem  \ref{thm2} is studied.

\begin{prop}\label{A}
If $\bu\in C^2(\bar{\Omega})$ satisfying $\inf_\Omega u_1>0$ and $\bu\cdot \bn=0$ on $\Gamma$, the Euler
system (\ref{Massconservation})-(\ref{Momentumconservation}) is
equivalent to the following system of equations
\begin{eqnarray}
\Div\,\ (\rho\mathbf{u})&=& 0,\label{0.1}\\
(\bu\cdot\nabla)B&=& 0,\label{0.2}\\
(\bu\cdot\nabla)\mathbf{\omega}+\mathbf{\omega}\Div\mathbf{u}-(\mathbf{\omega}\cdot\nabla)\bu&=&0,\label{0.3}
\end{eqnarray}
with boundary condition $\bu\cdot \bn =0$ on $\Gamma$ and
 the compatibility condition
\begin{equation}\label{0.4}
(\mathbf{u}\times
\omega)\cdot \tau=\nabla B\cdot \tau \,\
\text{on}\,\ \Gamma_-,
\end{equation}
where
$\mathbf{\omega}=\nabla\times\mathbf{u}$, $B$, and $\tau$ are the vorticity of the flow,  Bernoulli's function, and any tangential direction of $\Gamma_-$, respectively.
\end{prop}
\begin{pf}
The straightforward computation gives
\begin{equation}
(\mathbf{u}\cdot\nabla)\mathbf{u}=\nabla\left(\dfrac{1}{2}|\mathbf{u}|^2\right)-\mathbf{u}\times
(\nabla\times\mathbf{u}).\label{0}
\end{equation}
Hence
the momentum equations (\ref{Momentumconservation}) can be
written as
\begin{equation}
\mathbf{u}\times (\nabla\times\mathbf{u})=\nabla B.\label{1}
\end{equation}
This implies \eqref{0.2}.
Taking the operator $\nabla\times$ for (\ref{1}) yields
\begin{equation}
\nabla\times(\mathbf{u}\times (\nabla\times\mathbf{u}))=0.\label{3}
\end{equation}
With the aid of the identity
\begin{equation}
\nabla\times(\mathbf{u}\times\mathbf{\omega})=\mathbf{u}\Div\mathbf{\omega}
+(\mathbf{\omega}\cdot\nabla)\mathbf{u}-\mathbf{\omega}\Div\mathbf{u}
-(\mathbf{u}\cdot\nabla)\mathbf{\omega},\label{transport}
\end{equation}
and noting that $\text{div}(\nabla\times \bu)=0$,
one concludes that (\ref{3}) is equivalent to
(\ref{0.3}). The compatibility condition (\ref{0.4}) follows from (\ref{1}).

On the other hand, if $(\rho,\mathbf{u})$ is a solution of
(\ref{0.1})-(\ref{0.3}),  then it follows from (\ref{transport}) that \eqref{3} holds.
Therefore, there exists a smooth function $\Psi$ such that
\begin{equation}
\bu\times(\nabla\times\bu)=\nabla \Psi.\label{5}
\end{equation}
Obviously, $\bu\cdot \nabla \Psi=0$.
This, together with (\ref{0.2}), implies that
\begin{equation}
\bu\cdot\nabla (B-\Psi)\equiv 0.\label{6}
\end{equation}
Note that the compatibility condition (\ref{0.4}) yields
\begin{equation*}
\tau\cdot \nabla(B-\Psi)=0\,\ \text{on}\,\ \Gamma_{-}
\end{equation*}
for any tangential direction $\tau$ of $\Gamma_-$.
After choosing the constant suitably, one has $\Psi=B$ on $\Gamma_-$.
Since it is assumed that $\inf_\Omega u_1>0$, every point in $\Omega$ can be
traced back to a point on $\Gamma_{-}$ by a streamline. This gives $\Psi\equiv B$ in $\Omega$ and thus (\ref{1}).
Combining \eqref{1} with \eqref{0} gives the momentum equations \eqref{Momentumconservation}.
\end{pf}

%Noting that
%\[
%\bu \times \omega =(\bu_T+\bu_n)\times(\omega_T+\omega_n)=\bu_T\times\omega_T+\bu_T \times\omega_n +\bu_n\times \omega_T.
%\]
%Since $\bu_T\times \omega_n$ and $\bu_n\times \omega_T$ must be tangential vectors on the entrance,
%\[
%(\bu \times \omega)_T =\bu_T \times\omega_n +\bu_n\times \omega_T.
%\]
%If $B\equiv \bar B$ on $\Gamma_-$, then $\nabla_T B=0$ on $\Gamma_-$. Let $\omega_n=\tilde\kappa \bu_n$. Then $\omega_T=\tilde\kappa \bu_T$. Therefore, $\omega =\tilde\kappa u$ on $\Gamma_-$.

The following result on a div-curl system is used to solve the Euler system in Section \ref{Existence}.
\begin{lem}\label{5.2lemma2}
Let $\bw\in C^{1, \alpha}(\bar\Omega)$ and $\lambda \in C^{2,\alpha}(\bar{\Omega})$ satisfy $\min_{\bar{\Omega}} \lambda>0$, $\Div\,\ \bw=0$ in $\Omega$, and
\begin{equation}\label{5.2dc1}
\bw\times\mathbf{n}=0,\quad \mathbf{n}\cdot{\nabla  \lambda}=0\quad \text{on}\,\, \Gamma,
\end{equation}
where $\mathbf{n}$ is the outer unit normal of $\Gamma$.
If the vector field $\mathbf{v}=(v_1, v_2, v_3)\in C^{2, \alpha}(\bar \Omega)$ satisfies the compatibility conditions
\begin{equation}\label{5.2dc2}
v_i=\partial_{ii}v_i=0 \,\, \text{and}\,\, (\nabla\times\mathbf{v})\times\mathbf{e}_i =\mathbf{0} \quad \text{on}~\Gamma_i, \quad i=2, 3,
\end{equation}
where $\{\mathbf{e}_i\}_{i=1}^3$ is the canonical base for $\mathbb{R}^3$,
then the following div-curl system
\begin{equation}
\left\{
\begin{array}{llll}
\nabla\cdot(\lambda \mathbf{u})&= \nabla\cdot\mathbf{v}\,\,\ &\text{in}\,\ \Omega,\\
\nabla\times\mathbf{u}&=\bw \,\ &\text{in}\,\ \Omega,\\
\lambda \mathbf{u}\cdot\mathbf{n}&=\mathbf{v}\cdot\mathbf{n}\,\ &\text{on}\,\ \partial\Omega,
\end{array}\label{5.2B.1}
\right.
\end{equation}
admits a unique  solution $\mathbf{u}\in
C^{2,\alpha}(\bar{\Omega})$. Moreover, the solution $\mathbf{u}$ satisfies
\begin{equation}\label{estdc}
\|\mathbf{u}\|_{{k+1,\alpha};\bar{\Omega}}\leq C(\|\mathbf{v}\|_{{k+1,\alpha};\bar{\Omega}}+\|\bw\|_{{k,\alpha};\bar{\Omega}})\quad\text{for}\,\, k=0,1,
\end{equation}
where the constant $C$ depends only on $\|\lambda\|_{{k+1,\alpha};
{\Omega}}$ and $\min_{\bar{\Omega}}\lambda$.
\end{lem}
\begin{pf}
Step 1: Uniqueness. Suppose that $\mathbf{u}_i$ ($i=1,2$) are both solutions of the
problem (\ref{5.2B.1}). Set $\mathbf{U}=\mathbf{u}_1-\mathbf{u}_2$. Then
\begin{equation}
\left\{
\begin{array}{llll}
\nabla\cdot(\lambda \mathbf{U})=0\quad&\text{in}\quad\Omega,\\
\nabla\times\mathbf{U}=0\quad&\text{in}\quad\Omega,\\
\lambda \mathbf{U}\cdot\mathbf{n}=0\quad&\text{on}\quad\partial\Omega.
\end{array}\label{5.2B.2}
\right.
\end{equation}
It follows from the second equation in \eqref{5.2B.2} that there exists a function $\Phi$ such that $\mathbf{U}=\nabla\Phi$.
Furthermore, $\Phi$ solves the following co-normal boundary value problem for the elliptic equation,
\begin{equation*}
\left\{
\begin{array}{lll}
\nabla\cdot(\lambda \nabla\Phi)=0\quad&\text{in}\quad\Omega,\\
\lambda \dfrac{\partial\Phi}{\partial\mathbf{n}}=0\quad&\text{on}\quad
\partial\Omega.
\end{array}
\right.
\end{equation*}
Multiplying $\Phi$ on the both sides of the equation above and integrating by part yield
$$\int_{\Omega} \lambda |\nabla\Phi|^2=0.$$
Thus $\mathbf{U}=\nabla\Phi\equiv0$.

Step 2: Existence. First,  it follows from the notes in p.215 in \cite{GT} that the problem
\begin{equation}
\left\{
\begin{array}{lll}
\nabla\cdot(\lambda\nabla\phi)=\Div \mathbf{v}\quad&\text{in}\,\ \Omega,\\
\lambda \dfrac{\partial\phi}{\partial\mathbf{n}}=\mathbf{v}\cdot\mathbf{n}\quad&\text{on}\,\ \partial\Omega
\end{array}\label{5.2B.3}
\right.
\end{equation}
 has a unique solution $\phi\in H^1(\Omega)\cap\mcA$.
Furthermore, the Schauder estimate \cite[Theorem 6.2]{GT} gives that $\phi \in \mcC^{2, \alpha}(\Omega)$.
Since $\bv$ and $\lambda$ satisfy the compatibility conditions (\ref{5.2dc1})-(\ref{5.2dc2}),
 it follows from Lemma \ref{appendix1} in Appendix \ref{appendixreflection} that one has
$$\|\nabla\phi\|_{{k+1,\alpha};\bar{\Omega}}\leq
C\|\mathbf{v}\|_{{k+1,\alpha};\bar{\Omega}}$$
the constant $C$ depends only on $\|\lambda\|_{{k+1,\alpha;}{\Omega}}$ and $\min_{\bar{\Omega}}\lambda$.

Second, since $\Div \bw=0$,  it follows from ~\cite{Saranen} that  the problem
\begin{equation}
\left\{
\begin{array}{lll}
\nabla\times(\lambda^{-1}\nabla\times\mathbf{q})=\bw\quad&\text{in}\,\ \Omega,\\
\mathbf{n}\times\mathbf{q}=\mathbf{0}\quad&\text{on}\,\ \partial\Omega,\\
\nabla\cdot\mathbf{q}=0\quad&\text{on}\,\ \partial\Omega
\end{array}\label{5.2B.4}
\right.
\end{equation}
admits a unique weak solution  $q\in H^1(\Omega)$.
Furthermore, $\|\mathbf{q}\|_{H^1}\leq
C\|W\|_{L^2}$.

Note that the first equation in (\ref{5.2B.4}) can be rewritten as
\begin{equation}
-\lambda^{-1}\Delta\mathbf{q}+\nabla(\lambda^{-1})\times(\nabla\times\mathbf{q})=\bw,\label{5.2B.5}
\end{equation}
where we use identity $\nabla\times(\nabla\times\mathbf{q})\equiv\nabla\Div\mathbf{q}-\Delta\mathbf{q}$. Then it follows from \cite[Theorem 8.10]{GT} that $\mathbf{q}\in H^3(K)$ for any $K\Subset \Omega$.
Furthermore, the Schauder estimates for single
elliptic equation (\cite[Lemma 6.4]{GT}) yields $\mathbf{q}\in
C^{2,\alpha}(\bar{\Omega}\setminus (\mathcal{E}_\Gamma\cup \mathcal{E}_e))$. Furthermore,
\begin{equation*}
\|\mathbf{q}\|_{{k+2,\alpha};(K)}\leq
C_K(\|W\|_{{k,\alpha};\bar{\Omega}}+\|\mathbf{q}\|_{H^1(\Omega)})\,\
\text{for}\,\ K\Subset\bar{\Omega}\setminus (\mathcal{E}_\Gamma\cup \mathcal{E}_e).
\end{equation*}
Given a point $x\in (\mathcal{E}_\Gamma\cup \mathcal{E}_e)\setminus \mathscr{C}$,  without loss of generality, we choose $x=(0,1/2,0)$. Let us study the regularity of $\mathbf{q}$ near the point $x$. Define the odd and even extension $b\rightarrow \acute{b}$, $b\rightarrow \grave{b}$ respectively by
\begin{equation}\label{5.2extension}
\acute{b}(\cdot,x_3)=\text{sgn}(x_3)b(\cdot,|x_3|),\quad\grave{b}(\cdot,x_3)=b(\cdot,|x_3|).
\end{equation}
Now in $\tilde{\Omega}=[0,L]\times[0,1]\times[-1,1]$, we take $\tilde{\mathbf{q}}=(\acute{q}_1,\acute{q}_2,\grave{q}_3)$, $\tilde{W}=(\acute{W}_1,\acute{W}_2,\grave{W}_3)$. It follows from \eqref{5.2dc1} and the boundary condition for $\mathbf{q}$ in \eqref{5.2B.4} that
$\grave{\lambda}\in C^{2,\alpha}(\tilde{\Omega})$, $\tilde{W}\in C^{1,\alpha}(\tilde{\Omega})$ and $\tilde{\mathbf{q}}\in H^1(\tilde{\Omega})$. Furthermore, $\tilde{\mathbf{q}}$ solves the following problem in $\tilde{\Omega}$,
\begin{equation*}
\left\{
\begin{array}{lll}
\nabla\times(\grave{\lambda}^{-1}\nabla\times\tilde{\mathbf{q}})=\tilde{W}\quad\text{in}\quad\tilde{\Omega},\\
\nabla\cdot\tilde{\mathbf{q}} =0\quad\text{on}\quad\partial\tilde{\Omega},\\
\tilde{\mathbf{q}}\times\mathbf{n}=\mathbf{0}\quad\text{on}\quad\partial\tilde{\Omega}.
\end{array}
\right.
\end{equation*}
Note that $x$ is a regular boundary point on $\partial\tilde{\Omega}$, and thus applying the Schauder estimates near the boundary \cite[Lemma 6.4]{GT} gives
$$\|\tilde{\mathbf{q}}\|_{{k+2,\alpha};\bar{\tilde{\Omega}}\cap B_{1/3}(x)}\leq C(\|\tilde{\mathbf{q}}\|_{H^1(\tilde{\Omega})}+\|\tilde{\mathbf{w}}\|_{{k,\alpha};\bar{\tilde{\Omega}}})\quad \text{for}\,\, k=0,1.$$
Therefore, one has
\begin{equation}\label{estkey}
\|\mathbf{q}\|_{{k+2,\alpha};\bar{\Omega}\cap B_{1/3}(x)}\leq C(\|\mathbf{q}\|_{H^1(\Omega)}+\|\bw\|_{{k,\alpha};\bar{\Omega}})\quad \text{for}\,\, k=0,1.
\end{equation}

If $x\in \mathscr{C}$, after one more extension, $x$ becomes a regular boundary point and we can get an estimate similar to \eqref{estkey}. So one can conclude that $\mathbf{q}\in C^{k+2,\alpha}(\bar{\Omega})$.

 Let $\mathbf{r}=\lambda^{-1}\nabla\times\mathbf{q}$. Since $\mathbf{q}\times \mathbf{n}=\mathbf{0}$ on $\Gamma$, one has $\mathbf{q}\cdot \tau=0$ for any tangential direction $\tau$ on $\Gamma$. Given any closed set $S\subset \Gamma$, it follows from the Stokes theorem
 \[
 \int_S(\nabla\times\mathbf{q})\cdot\mathbf{n}d\sigma= \oint_{\partial S}\mathbf{q}\cdot\mathbf{\tau}ds=0
  \]
 that  $\mathbf{r}\cdot\mathbf{n}=0$. Hence $\mathbf{u}=\mathbf{v}+\mathbf{r}=\nabla\phi+\lambda^{-1}\nabla\times\mathbf{q}\in C^{k+1,\alpha}(\bar{\Omega})$ solves the problem \eqref{5.2B.1} and satisfies the estimate \eqref{estdc}.
\end{pf}

\section{Three Dimensional Potential Flows}\label{PotentialFlows}
Note that the compatibility condition \eqref{0.4} gives
\begin{equation}\label{vorticityIC}
\omega_2(0,x_2,x_3)=(-\partial_3B_0+u_2\kappa)/u_1,\quad\omega_3(0,x_2,x_3)=(\partial_2B_0+u_3\kappa)/u_1.
\end{equation}
If $B_0=\bar B$(constant) and $\inf_{\Omega}u_1>0$, the equation \eqref{0.2} shows  that $B\equiv\bar{B}$ in $\Omega$. Thus $\rho=H(\bar B-\frac{|\bu|^2}{2})$ where $H$ is the inverse function of $h$.
Therefore, if $\kappa\equiv0$ and $B_0\equiv\bar{B}$ on $\Gamma_-$, then one has $\omega(0,x_2,x_3)=0$. It follows from the transport equation \eqref{0.3} that $\omega\equiv0$ in $\Omega$. Therefore, there exists a function $\phi$ such that $\bu =\nabla\phi$. In order to fix the integral constant,  we always choose $\phi\in \mcA$. The continuity equation becomes
\begin{equation}
\nabla\cdot(\rho(|\nabla\phi|^2)\nabla\phi) =0,
\end{equation}
where $\rho(|\nabla\phi|^2) =H(\bar B-\frac{|\nabla\phi|^2}{2})$.

It is easy to compute that there exists a constant $c^*$ such that if $|\nabla\phi|< (>, =)c^*$, then $c(\rho)>(<, =)c^*$ (cf. \cite{CF}).  After nondimensionalization, we assume $c^*=1$. Thus the flow is subsonic (supersonic, sonic) for $|\nabla\phi|<1$( $|\nabla\phi|>1$, $|\nabla\phi|=1$) respectively.

We begin with the study on the wellposedness of the potential flows in bounded nozzles.
\begin{thm}\label{PotentialFlows0}
For any given $f$ satisfying \eqref{compatibility1} and \eqref{comp2}, there exists a
$\theta^*\in (0,\infty)$ such that for any $\theta\in (0, \theta^*)$, there exists a unique uniformly subsonic potential flow with boundary condition \eqref{Boundary2}. More precisely, there exists a unique solution $\phi \in C^{2, \alpha}(\bar\Omega)$ satisfying $\max_{\bar\Omega}|\nabla\phi|<1$ and
\begin{equation}\label{eq32}
\left\{
\begin{aligned}
&\nabla\cdot(\rho(|\nabla\phi|^2)\nabla\phi)=0,\,\quad \text{in}\,\, \Omega,\\
&\rho(|\nabla\phi|^2)\nabla\phi\cdot\bn=\theta f,\,\quad \text{on}\,\, \partial\Omega.
\end{aligned}
\right.
\end{equation}
Moreover, $\lim_{\theta\rightarrow \theta^*}\max_{\bar\Omega}|\nabla\phi|=1$. Finally, as $\theta\to \theta^*$, there exists a subsequence $\{\theta_n\}$ such that the corresponding solutions $\{(\rho_n, \bu_n)\}$ tend to $(\rho, \bu)$ which satisfies the Euler equations
\begin{equation}\label{wEulereq}
\left\{
\begin{aligned}
&\nabla\cdot(\rho(|\bu|^2)\bu)=0,\\
&\nabla\times\bu =0
\end{aligned}
\right.
\end{equation}
in the sense of distribution and the boundary condition
\begin{equation*}
\rho \bu \cdot \bn =\theta^*f
\end{equation*}
 in the sense of boundary trace.
\end{thm}

The main difficulties for the proof of Theorem \ref{PotentialFlows0} are that the equation becomes degenerate elliptic as $|\nabla\phi|\to 1$ and the boundary conditions are nonlinear.

Let $m\in \mathbb{N}$ and $\zeta_m$ be a smooth increasing function satisfying
\begin{equation}
\zeta_m(s)=\left\{
\begin{aligned}
&s,\,\, \text{if}\,\, s\leq 1-\frac{1}{m},\\
&1-\frac{2}{3m},\,\, \text{if}\,\, s\geq 1-\frac{1}{2m}.
\end{aligned}
\right.
\end{equation}

First, we have the following lemma on the truncated problem.
\begin{lem}\label{prioriestimates}
Suppose that $\phi\in \mcC^{2,\alpha}(\bar{\Omega})$ is a solution of the oblique problem
\begin{equation}\label{eq34}
\left\{
\begin{aligned}
&\nabla\cdot(\rho_m(|\nabla \phi|^2)\nabla\phi)=0\,\quad \text{in}\,\ \Omega,\\
&\rho_m(|\nabla \phi|^2)\nabla\phi\cdot \bn -\sigma \theta f(x)=0\,\quad  \text{on}\,\ \partial\Omega,
\end{aligned}
\right.
\end{equation}
where $\rho_m(s) =\rho(\zeta_m(s))$ and $\sigma\in [0, 1]$.
Then
\begin{equation}\label{est35}
\|\phi\|_{{2,\alpha};\bar{\Omega}}\leq C,
\end{equation}
where $C$ depends on $m$,  $\|f\|_{{1, \alpha};\partial\Omega}$ and monotonically on $\theta$ satisfying $C=0$ if $\theta=0$.
\end{lem}
\begin{pf}
The key step for the proof is the $L^\infty$ estimate for $\phi$.

Multiplying the both sides of the equation in \eqref{eq34} by $\phi$ and integrating by parts yield
\begin{equation}
\int_{\Omega} \rho_m |\nabla\phi|^2 dx = \int_{\partial \Omega} \sigma \theta f\phi dS \leq C \theta \|f\|_{L^2(\partial \Omega)}\|\phi\|_{L^2(\partial \Omega)}\leq C\theta\|f\|_{0;\partial\Omega}\|\nabla\phi\|_{L^2(\Omega)},
\end{equation}
where we use trace theorem and Poincar\'{e} inequality in the last inequality. Using Young inequality gives
$$\|\nabla\phi\|_{L^2(\Omega)}\leq C\theta\|f\|_{0;\partial\Omega}.$$
Thus it follows from Poincar\'{e} inequality that one has
\begin{equation}
\|\phi\|_{H^1(\Omega)}\leq C\theta \|f\|_{0;\partial\Omega}.
\end{equation}

Now we estimate $\|\phi\|_{L^{\infty}(\Omega)}$ by Moser's iteration.
Define $\hat{\phi}=\max\{\phi,1\}$. Multiplying the both sides of the equation in \eqref{eq34} by  $\psi=\hat{\phi}^q$ ($q\geq2$)
yields
\begin{equation}\label{eq38}
\int_{\Omega}q\rho_m\hat{\phi}^{q-1}|\nabla\hat{\phi}|^2 dx =\int_{\partial\Omega}\sigma \theta f\hat{\phi}^q dS \leq C \|f\|_{0;\partial\Omega} \|\hat\phi^q\|_{L^1(\partial \Omega)}.
\end{equation}
Note that the trace theorem gives $\|\psi\|_{L^1(\partial\Omega)}\leq C\|\psi\|_{W^{1,1}(\Omega)}$. Thus
\begin{equation}\label{eq39}
\|\hat\phi^q\|_{L^1(\partial \Omega)} \leq C\left(\int_{\Omega}\hat{\phi}^q dx+\int_{\Omega}q\hat{\phi}^{q-1}|\nabla\hat{\phi}|dx\right).
\end{equation}
Combining \eqref{eq38} and \eqref{eq39} together yields
\begin{equation}
\int_{\Omega}q\rho_m\hat{\phi}^{q-1}|\nabla\hat{\phi}|^2 dx
\leq C\|f\|_{0;\partial\Omega}\left(\int_{\Omega}q\hat{\phi}^q dx+\int_{\Omega}\epsilon q\hat{\phi}^{q-1}|\nabla\hat{\phi}|^2 dx \right),
\end{equation}
where we use the property $\hat\phi\geq 1$.
Thus
\begin{equation}
\int_{\Omega}|\nabla\hat{\phi}^{\frac{q+1}{2}}|^2 dx\leq C (q+1)^2\int_{\Omega}\hat{\phi}^{q+1}dx.
\end{equation}
Using Sobolev imbedding gives
\begin{equation}
\|\hat{\phi}\|_{L^{3(q+1)}}\leq \Bigl(C\|\hat{\phi}^{\frac{q+1}{2}}\|_{H^{1}(\Omega)}\Bigr)^{2/(q+1)} \leq  C^{\frac{1}{q+1}}(q+1)^{\frac{2}{q+1}}\|\hat{\phi}\|_{L^{q+1}}.
\end{equation}
When setting $q+1=3^{\nu}, \nu=1,2,3\cdots$ and letting $\nu\rightarrow\infty$, we have
$$\|\hat{\phi}\|_{L^{\infty}}=\lim_{\nu\to \infty} \|\hat{\phi}\|_{L^{3^{\nu}}}\leq C\|\hat{\phi}\|_{L^3}\leq C(\|\hat{\phi}\|_{L^{\infty}})^{\frac{1}{3}}\|\hat{\phi}\|_{L^2}^{\frac{2}{3}}.$$
Thus
 $$\|\hat{\phi}\|_{L^{\infty}}\leq C\|\hat\phi\|_{L^2}\leq C.$$
 Similarly, we have $\inf\check{\phi}\geq -C$ for $\check{\phi}=\min\{\phi,-1\}$. Therefore,
\[
\|\phi\|_{0;\bar\Omega}\leq C.
\]

Since $\partial_2\phi(x_1, 0, x_3)=0$ and $\partial_3\phi(x_1, x_2,0)=0$, we can extend $\phi$ to $[0, L]\times [-1,1]^2$ as follows
\begin{equation}\label{extension1}
\phi(x_1, x_2, -x_3)=\phi(x_1, x_2, x_3)\quad \text{and}\quad  \phi(x_1, -x_2, x_3) =\phi(x_1, x_2, x_3).
\end{equation}
Furthermore, we extend $f$ at the entrance
\begin{equation}
f(0, x_2, -x_3) =f(0, x_2, x_3),\quad\text{and}\quad  f(0, -x_2, x_3) =f(0, x_2, x_3)
\end{equation}
and the exit
\begin{equation}
f(L, x_2, -x_3) =f(L, x_2, x_3),\quad\text{and}\quad f(L, -x_2, x_3) =f(L, x_2, x_3),
\end{equation}
respectively. It is easy to see that periodic extension of $\phi$ with respect to $x_2$ and $x_3$ satisfies the same equation with periodic extension of $f$ in the domain $[0, L]\times \mathbb{R}^2$. After these extensions, there is no boundary in $x_2$ and $x_3$ direction.
The domain becomes a smooth domain.

As long as we have the $L^\infty$ estimate, it is the same to the estimate obtained in  \cite{LT} for nonlinear oblique derivative problems for quasilinear elliptic equations that  one has the estimate (\ref{est35}).
\end{pf}

 Set
\begin{equation}
E=\{\phi \in \mcC^{2, \alpha}(\bar\Omega)| F_m[\phi]=0, G_m[\sigma \theta, \phi]=0,\,\ \text{for some}\,\ \sigma\in[0,1]\},
\end{equation}
where
\[
F_m[\phi]=\nabla\cdot(\rho_m(|\nabla \phi|^2)\nabla\phi)\quad\text{and}\quad G_m[s, \phi]=\rho_m(|\nabla \phi|^2)\nabla\phi\cdot \bn -s f(x).
\]
It follows from Lemma \ref{prioriestimates} that $E$ is bounded in $C^{2,\alpha}(\bar{\Omega})$. By the method of continuity, the truncated problem \eqref{eq34} has a solution in $\mcC^{2, \alpha}(\bar\Omega)$.

To finish the proof of Theorem \ref{PotentialFlows0}, we first prove the uniqueness of solution for the problem \eqref{eq34}.  Suppose that there exist two uniformly subsonic  solutions $\phi_1$ and $\phi_2$. Set $\Phi=\phi_1-\phi_2$. We have
\begin{equation}\label{Eq}
\left\{
\begin{aligned}
& \partial_i (a_{ij}\partial_j \Phi)=0, \quad \text{in}\quad \Omega,\\
& a_{ij}\partial_j\Phi n_i =0, \quad \text{on}\quad \partial\Omega,
\end{aligned}
\right.
\end{equation}
where
\begin{equation}
a_{ij} =\int_0^1 [\rho(|\nabla (\phi_1+s \Phi)|^2)\delta_{ij} +2\rho'(|\nabla (\phi_1+s \Phi)|^2) \partial_i (\phi_1+s \Phi) \partial_j (\phi_1+s \Phi)] ds.
\end{equation}
Multiplying the both sides of the equation in \eqref{Eq} with $\Phi$ and integrating by parts yield
\begin{equation}
\int_{\Omega}|\nabla \Phi|^2ds \leq C \int_\Omega a_{ij}\partial_j\Phi\partial_i \Phi dx=0.
\end{equation}
Therefore, $\Phi\equiv C$. Recall that $\Phi\in \mcA$, thus $\Phi=0$. This gives the uniqueness of the solution.

Let $\phi_m(\cdot, \theta)$ be the unique solution of the problem
\[
F_m[\phi]=0\quad \text{in}\quad \Omega \quad \text{and}\quad G_m[\theta, \phi]=0\quad \text{on}\quad \partial\Omega.
\]
Define $M_m(\theta)=\sup_{\bar{\Omega}}|\nabla\phi_m(\cdot;\theta)|^2$.
Since $\phi_m\in C^{2,\alpha}(\Omega)$, it follows from the Schauder estimate \cite{GT}  that $|M_m(\theta)-M_m(\tilde\theta)|\leq C_m(|\theta-\tilde\theta|)$. Thus, $M_m(\theta)$ is a continuous function of $\theta$. Set $\theta_m=\sup\{s\in[0,1]| M_m(\tau)\leq 1-\frac{1}{m}\,\ \text{for all}\,\ \tau\in [0,s)\}$. Note that
$M_m(0)=0$. Thus, $\theta_m$ is well defined.
It is obvious that $\{\theta_m\}$ is an increasing sequence and $M_m(\theta_m)=1-\frac{1}{m}$ due to the continuity of $M_m$. Take $\theta^*=\lim_{m\rightarrow\infty}\theta_m$, and then $\theta^*$ is the exact critical number stated in Theorem \ref{PotentialFlows0}.

The existence of weak solution for \eqref{wEulereq} follows from the compensated compactness argument in \cite{Huang}.

Thus, the proof of Theorem \ref{PotentialFlows0} is completed.

%\newpage
In fact, under the assumption in Theorem \ref{PotentialFlows0}, we can get better regularity for the solution. Furthermore, in order to show that the potential model is equivalent to the original Euler flows, we need to prove the positivity of the velocity component $u_1$. These are stated in the following proposition.
\begin{prop}\label{prop8}
The solution obtained in Theorem \ref{PotentialFlows0} satisfies $\phi\in C^{3, \alpha}(\bar\Omega)$ and $u_1>0$.
\end{prop}
\begin{pf}
Taking derivative for the equation in \eqref{eq32} with respect to $x_2$ yields that $\partial_2\phi$ satisfies the equation
\begin{equation}\label{eq47}
\partial_i((\rho \delta_{ij}+2\rho'\partial_i\phi\partial_j\phi)\partial_j(\partial_2\phi))=0.
\end{equation}
Differentiating the boundary condition in \eqref{eq32} on the boundary $\Gamma_{\pm}$ with respect to $x_2$ gives
\begin{equation}\label{eq48}
(\rho+2\rho'(\partial_2\phi)^2)\partial_1\partial_2\phi +\sum_{i=2,3}2\rho'\partial_1\phi\partial_i\phi\partial_i\partial_2\phi=\theta\partial_2 f.
\end{equation}
Similarly, $\partial_2\phi$ satisfies the similar boundary condition on $\Gamma_3$.
Note that $\partial_2\phi=0$ on $\Gamma_2$. Under the extension \eqref{extension1}, $\partial_2\phi$ is a solution of the problem  \eqref{eq47} and \eqref{eq48} in the region $[0, L]\times [-1, 1]\times [0,1]$. It is easy to check that \eqref{eq48} is an oblique derivative boundary condition for $\partial_2\phi$. The Schauder estimate for elliptic equation with oblique derivative boundary conditions \cite{GT} gives
\begin{equation}\label{eq49}
\|\partial_2\phi\|_{{2,\alpha};\Omega}\leq C\|f\|_{{2,\alpha};\partial\Omega}.
\end{equation}
Similarly, one has
\begin{equation}\label{eq50}
\|\partial_3\phi\|_{{2,\alpha};\Omega}\leq C\|f\|_{{2,\alpha};\partial\Omega}
\end{equation}
Obviously, $\partial_1\phi$ also satisfies the equation \eqref{eq47}. Now let us derive the boundary condition for $\partial_1\phi$ on $\Gamma_-$.
It follows from the equation \eqref{eq32} that $\partial_1\phi$ satisfies
\begin{equation}\label{Add}
\begin{aligned}
&[\rho+2\rho'(\partial_1\phi)^2]\partial_1\partial_1\phi +\sum_{i=2}^3 2\rho'\partial_1\phi\partial_i\phi\partial_i\partial_1\phi = -\sum_{i, j=2}^3 [\rho\delta_{ij}+2\rho'\partial_i \phi\partial_j\phi]\partial_{ij}\phi.
\end{aligned}
\end{equation}
This is also an oblique derivative boundary condition for $\partial_1\phi$. Similarly, $\partial_1\phi$ satisfies the same boundary condition (\ref{Add}) on $\Gamma_+$. On $\Gamma_i$, it holds that $\partial_i\partial_1\phi$=0. Thus, the Schauder estimate yields
\begin{equation}\label{eq52}
\|\partial_1\phi\|_{{2,\alpha};\Omega}\leq C \sum_{i=2}^3\|\partial_i \phi\|_{{2,\alpha};\Omega}.
\end{equation}
Combining the estimate \eqref{eq49}-\eqref{eq52} together, we have
\begin{equation}
\|\nabla\phi\|_{{2, \alpha};\bar\Omega}\leq C \|f\|_{{2, \alpha};\partial\Omega}.
\end{equation}
This, together with Theorem \ref{PotentialFlows0}, implies
\begin{equation}
\|\phi\|_{{3, \alpha};\bar\Omega}\leq C \|f\|_{{2, \alpha};\partial\Omega}.
\end{equation}
Thus $\phi\in C^{3, \alpha}(\bar\Omega)$.

Now let us prove the positivity of velocity component $u_1=\partial_1\phi$.
Straightforward computations show that $u_1$ is a solution of the following problem
\begin{equation}
\left\{
\begin{aligned}
&\partial_i((\rho(|\nabla\phi|^2)\delta_{ij}+2\rho'(|\nabla\phi|^2)\partial_i\phi\partial_j\phi)\partial_j u_1)=0,\,\quad \text{in}\,\ \Omega,\\
&\nabla u_1\cdot\bn=0,\,\ \text{on}\,\,\quad \Gamma,\\
&u_1=\frac{\pm \theta f}{\rho},\,\ \text{on}\,\,\quad \Gamma_{\pm}.
\end{aligned}
\right.
\end{equation}
The maximum principle implies that $u_1$ achieves its minimum on the boundary.
By Hopf Lemma, the minimum of $u_1$ cannot be achieved on $\Gamma$. Meanwhile, on $\Gamma_{\pm}$, $u_1\geq C>0$ with the constant $C$ depending on $\inf_{\Gamma_-\cup\Gamma_+}|f|$. Therefore, $u_1>0$ on $\partial\Omega$.
This finishes the proof of the proposition.
\end{pf}

Theorem \ref{PotentialFlows0}, together with Proposition \ref{prop8}, gives Theorem \ref{thm1}.

%\newpage

\section{ General Three Dimensional Flows}\label{Existence}
In this section, we prove Theorem \ref{thm2}. By the assumption of the theorem, there exists a subsonic flow $(\bar{\rho},\bar{\mathbf{u}})$ with $\bar{\mathbf{u}}=\nabla\bar{\phi}$ and $h(\bar{\rho})+\frac{1}{2}|\bar{\mathbf{u}}|^2=\bar{B}$ satisfying
\begin{equation}\label{potential1}
\left\{
\begin{aligned}
&\nabla\cdot(\rho(|\nabla\bar{\phi}|^2)\nabla\bar{\phi})=0\quad\text{in}\quad\Omega,\\
&\rho(|\nabla\bar{\phi}|^2)\nabla\bar{\phi}\cdot\mathbf{n}=\bar{f}\quad\text{on}\quad\partial\Omega,
\end{aligned}
\right.
\end{equation}
where $\bar{f}$ satisfies the conditions (\ref{compatibility1}), (\ref{positive}) and (\ref{comp2}). According to the argument in Section \ref{PotentialFlows}, we have showed that
\begin{equation}\label{compatibility}
\bar{\phi}\in C^{3,\alpha}(\bar{\Omega})\quad \text{and}\quad  \inf_{\Omega}\bar{u}_1>0.
\end{equation}
Furthermore, we can extend $\bar\phi$ to be periodic with respect to $x_2$ and $x_3$ in the domain $([0,L]\times\mathbb{R}^2)$. Consequently,
\begin{equation}\label{pcomp}
\partial_{i}^3\bar{\phi}=0\quad\text{on}\quad\Gamma_i~(i=2,3).
\end{equation}
\begin{pf2}
Denote $\sigma_0=\inf_{\Omega}\bar{u}_1$, $\sigma_1=2\|\bar{\mathbf{u}}\|_{C^{2,\alpha}(\bar{\Omega})}$.
%First we modify the boundary data $h$ and $B_0$ by $h_{\mu}$ and
%$B_{0,\mu}$ such that
%\begin{equation}
%\begin{array}{lll}
%h_{\mu}\in C^{\infty}_{c}(\Gamma_-),&
%\|h_{\mu}\|_{C^{\alpha}(\Gamma_-)}\leq
%C_0\|h\|_{C^{\alpha}(\Gamma_-)},\\
%B_{0,\mu}\in C^{\infty}_{c}(\Gamma_-),& \|\nabla B
%_{0,\mu}\|_{C^{\alpha}(\Gamma_-)}\leq C_0\|\nabla B_0
%\|_{C^{\alpha}(\Gamma_-)}.
%\end{array}
%\end{equation}
%For simplicity, we require the compact supports for $h_{\mu}$ and
%$\nabla B_{0,\mu}$ are just a little larger than those for $h$ and
%$\nabla B_0$, i.e., $h_{\mu}\equiv0, \nabla
%B_{0,\mu}\equiv\b0\quad
%\text{for}\quad\gamma^2_{2}+\gamma^2_{3}\geq\delta_0$. We still
%denote the constant $C$ for it is independent of $\mu$.
Set
$$\mathcal{S}=\{\mathbf{u}\in C^{2,\alpha}\bar{\Omega}):
\|\mathbf{u}-\bar{\mathbf{u}}\|_{{2,\alpha;}\bar{\Omega}}\leq\sigma,
\mathbf{u}\cdot\mathbf{n}=0 \,\,\text{and}\,\,  (\nabla\times\mathbf{u})\times\mathbf{n}=0\,\ \text{on}\,\ \Gamma \}.$$
The constant $\sigma \leq\sigma_0/2$ is to be determined so that $u_1>\sigma_0/2$ and
$\|\mathbf{u}\|_{C^{2,\alpha}(\bar{\Omega})}\leq\sigma_0+\sigma_1/2\leq\sigma_1$ for any $\mathbf{u}\in\mathcal{S}$.
Given $\mathbf{u}\in\mathcal{S}$, we construct a vector $\mathbf{v}=\mathcal{T}\mathbf{u}$ by the following six steps.  The key idea is to show that $\mathcal{T}$ is a contract map from $\mathcal{S}$ to itself.

{Step 1:} Streamlines and solving Bernoulli function. First we extend the
function $\dfrac{u_i}{u_1} (i=2,3)$ by $U_i\in
C^{2,\alpha}([0,L]\times\mathbb{R}^2)$ with
$$\|U_i\|_{{2,\alpha;}[0,L]\times\mathbb{R}^2}\leq C_0\left\|\dfrac{u_{i}}{u_{1}}\right\|_{{2,\alpha;}\bar{\Omega}}$$
for some fixed constant $C_0$.
The streamline is defined by
\begin{equation}\label{characteristicline}
\left\{
\begin{aligned}
&\dfrac{dX_i}{ds}(s; x)=U_i(s,X_2(s; x),X_3(s; x)),\\
&X_i(s=x_1)=x_i.
\end{aligned}
\right.
\end{equation}
Later on, we denote the streamline to be $(x_1,X_2(s;x_1,x_2,x_3),X_3(s;x_1,x_2,x_3))$.
Thanks to the slip boundary condition on $\Gamma$, i.e., $\mathbf{u}\cdot\mathbf{n}=0$, the local uniqueness guarantees that
the particle path cannot intersect with the nozzle wall provided that $(x_1,x_2,x_3)\in\Omega$; while for any $(x_1,x_2,x_3)\in\Gamma$, the streamline always lies on the nozzle wall $\Gamma$. Hence, each streamline can
 be uniquely traced back to some point belonging to $\Gamma_-$. Let $\gamma_i(x_1,x_2,x_3)=X_i(0;x_1,x_2,x_3)$.
Then  we have the following estimates.
\begin{lem}\label{lem1}
(1) The streamline is also $C^{2,\alpha}$-smooth provided that the velocity field has $C^{2,\alpha}$ regularity, i.e.,
\begin{equation}\label{tranlabel1}
\sum_{i=2}^3\|X_i(s; \cdot)\|_{{2,\alpha};\bar{\Omega}}\leq C,
\end{equation}
with $C$ depending only on $\|\bu\|_{C^{2,\alpha}(\bar{\Omega})}$.

(2) Suppose $X_i(s)$ and $\tilde{X}_i(s)$ are streamlines associated with $\bu$ and $\tilde{\bu}$,
then
\begin{equation}\label{tranlabel2}
\sum^{3}_{i=2}\|X_i(s)-\tilde{X}_i(s)\|_{{1,\alpha};\bar{\Omega}}\leq C\|\bu-\tilde{\bu}\|_{{1,\alpha};\bar{\Omega}};
\end{equation}
where $C$ depends on $\|\bu\|_{C^{2,\alpha}(\bar{\Omega})}$ and $\|\tilde\bu\|_{C^{2,\alpha}(\bar{\Omega})}$.

(3) On the solid boundary $\Gamma$, we have
\begin{equation}\label{tranlabel3}
\partial_{x_i} X_j(s;x_1,x_2,x_3)\equiv 0\quad \text{for}\,\, x\in \Gamma_i, \,\, i, j=2, 3\,\, \text{and}\,\,i\neq j.
\end{equation}
\end{lem}
The proof for Lemma \ref{lem1} is given in the Appendix \ref{appendixstream}.

Using Lemma \ref{lem1} gives
\begin{equation}\label{characteristicline2}
\begin{aligned}
&\sum_{i=2}^{3}\|\gamma_i\|_{{2,\alpha};\bar{\Omega}}\leq C_1\quad\text{and}\,\,
\sum_{i=2}^{3}\|\gamma_i-\tilde{\gamma}_i\|_{{1,\alpha};\bar{\Omega}}\leq C_2\|\bu-\tilde{\bu}\|_{{1,\alpha};\bar{\Omega}}.
\end{aligned}
\end{equation}
where $C_1$ depends on $\|\bu\|_{{2,\alpha};\bar{\Omega}}$ and $C_2$ depends on $\|\bu\|_{{2,\alpha};\bar{\Omega}}$ and $\|\tilde{\bu}\|_{{2,\alpha};\bar{\Omega}}$. From now on, the constants $C$,  $C_1$ and $C_2$ may change from line to line. However, they keep the property that $C$ depends only on $\|\bar\bu\|_{2, \alpha;\bar \Omega}$, $C_1$ depends only on  $\|\bu\|_{{2,\alpha};\bar{\Omega}}$ and $C_2$ depends on $\|\bu\|_{{2,\alpha};\bar{\Omega}}$ and $\|\tilde{\bu}\|_{{2,\alpha};\bar{\Omega}}$.
Furthermore, it follows from from Lemma \ref{lem1} that
\begin{equation}\label{sdc}
\begin{array}{lll}
\partial_2\gamma_3(x_1,0,x_3)=\partial_2\gamma_3(x_1,1,x_3)=0 \,\,\text{and}\,\,
\partial_3\gamma_2(x_1,x_2,0)=\partial_3\gamma_2(x_1,x_2,1)=0.
\end{array}
\end{equation}

Note that $\dfrac{dB}{ds}(s; x)=0$, so the Bernoulli's function is invariant along the streamline. Therefore,
$$B(x_1,x_2,x_3)=B(s,X_2(s;x_1,x_2,x_3),X_3(s;x_1,x_2,x_3)).$$
In particular,
$$B(x_1,x_2,x_3)=B(0,X_2(0;x_1,x_2,x_3),X_3(0;x_1,x_2,x_3))=B_0(\gamma_2(x_1,x_2,x_3),\gamma_3(x_1,x_2,x_3)).$$
Since $B_0(0,0)=\bar{B}$, this, together with \eqref{characteristicline2}, yields
$$\|B-\bar{B}\|_{{2,\alpha};\bar{\Omega}}\leq C_1\|B_0-\bar{B}\|_{{2,\alpha};\Gamma_-}\leq C_1 \epsilon.$$
It follows from (\ref{sdc}) and the condition (\ref{B0}) that one has
\begin{equation}\label{Bcomp}
\mathbf{n}\cdot\nabla B=0\quad\text{on}\quad\Gamma.
\end{equation}
Let $\tilde{B}(x_1, x_2, x_3)$ be the Bernoulli function associated with $\tilde{\bu}$. Then we have
\begin{equation}\label{Bernoulli}
\begin{aligned}
\quad &\|B-\tilde{B}\|_{{1,\alpha};\bar{\Omega}}\leq C_1\|B_0-\bar B\|_{{2,\alpha};\Gamma_-}\Bigl(\sum^{3}_{i=2}\|\gamma_i-\tilde{\gamma}_i\|_{{1,\alpha};\bar{\Omega}})\Bigr)\\
\leq &C_1\epsilon\sum^{3}_{i=2}\|\gamma_i-\tilde{\gamma}_i\|_{{1,\alpha};\bar{\Omega}}\leq C_2\epsilon \|\bu-\tilde{\bu}\|_{{1,\alpha};\bar{\Omega}}.
\end{aligned}
\end{equation}

Step 2: Solving the vorticity. Extend the matrix
$\dfrac{(\Div\mathbf{u})I-(\nabla\mathbf{u})^{T}}{u_1}$ by $V\in
C^{1,\alpha}([0,1]\times\mathbb{R}^2)$. Let $\Lambda=(\Lambda_1,\Lambda_2,\Lambda_3)^T$ be the
solution to the following linear ODE system.
\begin{equation}\label{5.4transport}
\left\{
\begin{aligned}
&\dfrac{d\Lambda(s,X_2(s;x),X_3(s;x))}{ds}+V(s,X_2(s;x),X_3(s;x))\Lambda(s,X_2(s;x),X_3(s;x))=0,\\
&\Lambda(0,X_2(0;x),X_3(0;x))=\Bigl(-\kappa,\frac{\kappa u_2-\partial_3B_0}{u_1},\frac{\kappa u_3+\partial_2B_0}{u_1}\Bigr)(\gamma_2(x),\gamma_3(x))\,\
\text{at}\,\ \Gamma_-
\end{aligned}\right.
\end{equation}
for $x=(x_1,x_2,x_3)$.
Take $\omega(x_1,x_2,x_3)=\Lambda(x_1, X_2(x_1; x), X_3(x_1; x))$. Then
\begin{equation}\label{t1}
\omega(x)=e^{-\int^{x_1}_{0}V(\tau,X_2(\tau),X_3(\tau))d\tau}\Lambda_0,
\end{equation}
where $\Lambda_0= \Lambda(0,X_2(0;x),X_3(0;x))$.
This
yields that
\begin{equation}
\|\omega\|_{C^{1,\alpha}(\bar{\Omega})}\leq C_1(\|\nabla
B_0\|_{C^{1,\alpha}(\Gamma_-)}+\|\kappa\|_{C^{1,\alpha}(\Gamma_-)}).
\end{equation}

Note that $\omega$ defined in \eqref{t1} satisfies
\begin{equation}\label{tt1}
\mathbf{u}\cdot\nabla\omega+(\Div\mathbf{u})\omega-(\omega\cdot\nabla)\mathbf{u}=0.
\end{equation}
Since $\mathbf{u}\cdot\mathbf{n}=0$ and $(\nabla\times\mathbf{u})\times\mathbf{n}=0$ holds on $\Gamma$, one has $\partial_2u_j=\partial_ju_2\equiv0$  for $j=1,3$ on $\Gamma_2$. Thus, on $\Gamma_2$, the equation \eqref{tt1} can be written as
\begin{equation}\label{ttt1}
\mathbf{u}\cdot\nabla\omega_i+(\Div\mathbf{u})\omega_i- \sum_{j=1, 3} \omega_j \partial_j u_i=0\quad \text{for}\,\, i=1, 3.
\end{equation}
Since  $\kappa=\partial_2 B_0 =0$ on $\Gamma_-\cap{ \Gamma_2}$, one has $\omega_i\equiv0(i=1,3)$ on $\Gamma_-\cap{ \Gamma_2}$. Solving the equations \eqref{ttt1} on $\Gamma_2$ yields  that $\omega_i=0 (i=1,3)$ on $\Gamma_2$. Similarly, $\omega_1=\omega_2=0$ on $\Gamma_3$. Therefore, we have
\begin{equation}\label{0vorticity}
\omega\times\mathbf{n}=0\,\ \text{on}\,\ \Gamma.
\end{equation}

Furthermore, let $\omega$ and $\tilde{\omega}$ be the solutions defined in \eqref{t1} associated with $\bu$ and $\tilde{\bu}$, respectively. Then
\begin{equation}\label{vorticity}
\begin{aligned}
\|\omega-\tilde{\omega}\|_{{\alpha};\bar{\Omega}}\leq&
C_2(\|\nabla_{T}\kappa\|_{{\alpha};\Gamma_-}+\|\nabla_{T}^{2}B_0\|_{{\alpha};\Gamma_-})\|\bu-\tilde{\bu}\|_{{1,\alpha};\Gamma_-}\\
\leq&C_2\epsilon\|\bu-\tilde{\bu}\|_{{1,\alpha};\bar{\Omega}}.
\end{aligned}
\end{equation}

{Step 3}: $\omega$ is divergence free. Since $(\mathbf{u}\times\mathbf{\omega})\cdot\tau=\partial_{\tau}B_0$ for any tangential direction $\tau$ on $\Gamma_-$, the Stokes theorem gives
$$\int_S[\nabla\times(\mathbf{u}\times\mathbf{\omega})]\cdot d\mathbf{S}=\oint_{\partial S}(\mathbf{u}\times\mathbf{\omega})\cdot d\mathbf{l}=0.$$
This yields  that
$$\mathbf{n}\cdot \left(\nabla\times(\mathbf{u}\times\mathbf{\omega})\right)=0\quad \text{on}\,\, \Gamma_-.$$
Combining \eqref{tt1} with the identity \eqref{transport} gives $(\mathbf{n}\cdot\mathbf{u})\Div\mathbf{\omega}=0$ on $\Gamma_-$. Since $\mathbf{n}\cdot(\rho\mathbf{u})=f<0$ on $\Gamma_-$, one has
$\Div\mathbf{\omega}=0$ {on} $\Gamma_{-}$.

Differentiating (\ref{tt1}) yields
\begin{equation}\label{5.3.215}
0=\Div((\mathbf{u}\cdot\nabla)\mathbf{\omega}+\mathbf{\omega}\Div\mathbf{u}-(\mathbf{\omega}\cdot\nabla)\mathbf{u})
=(\mathbf{u}\cdot\nabla)\Div\mathbf{\omega}+(\Div\mathbf{u})(\Div\mathbf{\omega}).
\end{equation}
This, together with the boundary condition for $\Div \omega$ at $\Gamma_-$, gives $\Div\mathbf{\omega}\equiv 0$ in
$\Omega$.

{Step 4}: Construct the vortical part of the velocity field. It follows from the compatibility condition
\begin{equation}
\mathbf{u}\cdot\mathbf{n}=\nabla B\cdot\mathbf{n}=0,\quad (\nabla\times\mathbf{u})\times\mathbf{n}=\mathbf{0}\,\,\text{on}\,\, \Gamma
\end{equation}
that $\rho =H(B-\frac{1}{2}|\bu|^2)$ satisfies $\nabla\rho\cdot\mathbf{n}=0$ on $\Gamma$. Furthermore, we have  $\omega\times\mathbf{n}=0$ on $\Gamma$. Hence Lemma \ref{5.2lemma2} guarantees that there exists a unique  $W\in C^{2,\alpha}(\bar{\Omega})$ satisfying
\begin{equation}\label{eqW}
\left\{
\begin{aligned}
&\Div(\rho W)=0\quad\text{in}\,\ \Omega,\\
&\nabla\times W=\omega\quad\text{in}\,\ \Omega,\\
&\rho W\cdot\mathbf{n}=0\quad\text{on}\,\ \partial\Omega.
\end{aligned}
\right.
\end{equation}
Moreover, $W$ satisfies the following estimate
\begin{equation*}
\|W\|_{{2,\alpha};\bar{\Omega}}\leq C_1\|\omega\|_{{1,\alpha};\bar{\Omega}}\leq C_1\epsilon.
\end{equation*}
Let $W$ and $\tilde{W}$ be solutions associated with $\mathbf{u}$ and $\tilde{\mathbf{u}}$, respectively. Then $W-\tilde{W}$ solves the system
\begin{equation}\label{add0}
\left\{
\begin{aligned}
&\Div(\rho (W-\tilde{W}))=\Div((\tilde{\rho}-\rho)\tilde{W})\quad&&\text{in}\quad\Omega,\\
&\nabla\times (W-\tilde{W})=\omega-\tilde{\omega}\quad&&\text{in}\quad\Omega,\\
&\rho (W-\tilde{W})\cdot\mathbf{n}=(\tilde{\rho}-\rho)\tilde{W}\cdot\mathbf{n}\quad&&\text{on}\quad\partial\Omega.
\end{aligned}
\right.
\end{equation}

It is easy to see  that $\nabla \rho \cdot \mathbf{n}=0$ and $(\omega-\tilde\omega)\times\mathbf{n}=0$ on $\Gamma$. Differentiating \eqref{eqW} and the straightforward computations yield
\begin{equation*}
(\tilde\rho-\rho)\tilde W_i=\partial_{ii}[(\tilde\rho-\rho)\tilde W_i]=0\,\, \text{and}\,\, \{\nabla \times [(\tilde\rho-\rho)\tilde W]\}\times \mathbf{e}_i=0 \,\, \text{on}\,\, \Gamma_i,\,\, i=2, 3.
\end{equation*}
  Therefore, applying the priori estimate given in Lemma \ref{5.2lemma2} yields
\begin{equation}\label{5.4contraction3}
\|W-\tilde{W}\|_{{1,\alpha}}\leq C_1(\|\omega-\tilde{\omega}\|_{{\alpha}}+\|\tilde{W}\|_{{2,\alpha}}\|\rho-\tilde{\rho}\|_{{2,\alpha}})
\leq C_2\epsilon\|\mathbf{u}-\tilde{\mathbf{u}}\|_{{1,\alpha}}
\end{equation}

{Step 6}: Construct the irrotational part of the velocity field.

Note that the following compatibility conditions hold at $\Gamma$, i.e.,
\begin{equation}\label{5.4vorticity}
\nabla B\cdot\mathbf{n}=W\cdot\mathbf{n}=0,\quad(\nabla\times W)\times\mathbf{n}=0.
\end{equation}
We consider the
following nonlinear problem
\begin{equation}\label{5.4nonlinear}
\left\{
\begin{aligned}
&\Div(H(B-\dfrac{1}{2}|\nabla\phi+W|^2)\nabla\phi)=0\quad&&\text{in}\,\ \Omega,\\
&H(B-\dfrac{1}{2}|\nabla\phi+W|^2)\nabla\phi\cdot\mathbf{n}=f\quad&&\text{on}\,\
\partial\Omega.
\end{aligned}
\right.
\end{equation}
Then we have the following results on the problem \eqref{5.4nonlinear}.
\begin{lem}\label{5.4nonlinearlemma}
There exist two positive constants $\delta_1$ and $\sigma_2 \in (0, \sigma_0/2)$ such that if
$$\|f- \bar{f}\|_{{2,\alpha};\partial\Omega}\leq\delta_1, \quad \|B-\bar{B}\|_{{2,\alpha};\bar{\Omega}}+\|W\|_{{2,\alpha};\bar{\Omega}}\leq\sigma_2,$$
then there exists a solution $\phi\in \mcC^{3,\alpha}(\bar{\Omega})$ of (\ref{5.4nonlinear}) satisfying
\begin{equation}
\|\nabla(\phi-\bar{\phi})\|_{{2,\alpha};\bar{\Omega}}\leq
C_1(\|B-\bar B\|_{{2,\alpha};\bar{\Omega}}+ \|W\|_{{2,\alpha};\bar{\Omega}} +\|f-\bar{f}\|_{{2,\alpha};\partial\Omega}).
\end{equation}
Furthermore, suppose $\phi$ and $\tilde{\phi}$ are solutions associated with $(W,B)$ and $(\tilde{W},\tilde{B})$, respectively. Then, $\phi-\tilde{\phi}$ satisfies
the estimate
\begin{equation}\label{5.4contraction4}
\|\nabla(\phi-\tilde{\phi})\|_{{1,\alpha};\bar{\Omega}}\leq C_2(\|B-\tilde{B}\|_{{1,\alpha};\bar{\Omega}}+\|W-\tilde{W}\|_{{1,\alpha};\bar{\Omega}}).
\end{equation}
\end{lem}

We first use Lemma \ref{5.4nonlinearlemma} to finish the proof of Theorem \ref{thm2}. The proof of Lemma \ref{5.4nonlinearlemma} is given at the end of this section.

\text{Step 5}: Existence of the solution of the Euler system

Define
$\mathbf{v}=\mathcal{T}\mathbf{u}=\nabla\phi+W$. If $\epsilon$ is sufficiently small, then one has
$$\|\nabla(\phi-\bar{\phi})\|_{{2,\alpha};\bar{\Omega}}\leq\sigma_0/2\quad\text{and}\,\, \|W\|_{{2,\alpha};\bar{\Omega}}\leq\sigma_0/2.$$
Thus $\|\bv-\bar \bu\|_{2, \alpha;\Omega} \leq \sigma_0$. Note that
$$\nabla\phi\cdot\mathbf{n}=W\cdot\mathbf{n}=0\quad \text{on}\quad \Gamma$$
yields $\bv \cdot \mathbf{n}=0$. Furthermore, $\nabla \times\bv =\nabla \times W$ gives $(\nabla \times \bv)\times \mathbf{n}=0$ on $\Gamma$. Thus $\mathbf{v}\in\mathcal{S}$.
Furthermore, given $\bu^{(0)}\in \mathcal{S}$, we can define $\bu^{(n)}=\mathcal{T}\bu^{(n-1)}$ for $n\geq 1$. It is easy to see that $\bu^{(n)}\in \mathcal{S}$. Furthermore, combining the estimates \eqref{Bernoulli}, \eqref{0vorticity}, \eqref{5.4contraction3}, and \eqref{5.4contraction4} together and choosing $\epsilon$ sufficiently small yield
\begin{equation}
\|\bu^{(n)}-{\bu^{(n-1)}}\|_{1, \alpha; \Omega}\leq \frac{1}{2}\|\bu^{(n-1)} -{\bu^{(n-2)}}\|_{1, \alpha; \Omega}.
\end{equation}
Thus $\bu^{(n)}$ converges to $\bu$ in $C^{1, \alpha}(\bar\Omega)$. Since $\bu^{(n)}\in \mathcal{S}$, one has $\bu \in \mathcal{S}$. Therefore,
there exists a unique fixed point $\mathbf{u}$ for $\mathcal{T}$ in $\mathcal{S}$.
Furthermore, $(\mathbf{u},\omega,B)$
satisfies (\ref{0.1})-(\ref{0.3}). Thanks to Proposition \ref{A}, $\mathbf{u}$ is a solution of
the Euler system
(\ref{Massconservation})-(\ref{Momentumconservation}) with the
condition (\ref{Boundary}) and (\ref{Entrance2}). This finishes
the proof of Theorem \ref{thm2}.
\end{pf2}

In order to complete the proof of Theorem \ref{thm2}, we need only to prove Lemma \ref{5.4nonlinearlemma}.

Proof of Lemma \ref{5.4nonlinearlemma}.
Take $\Phi=\phi-\bar{\phi}$ where $\bar{\phi}$ is  the corresponding potential flow.
Then $\Phi$ satisfies the system
\begin{equation}\label{5.4psystem}
\left\{
\begin{aligned}
&\Div(\bar{\mathcal{M}}\nabla \Phi)=\Div\mathcal{F}(x,B,W,\nabla\Phi)\quad\text{in}\,\ \Omega,\\
&\bar{\mathcal{M}}\nabla \Phi\cdot\mathbf{n}=\mathcal{F}(x,B,W,\nabla\Phi)\cdot\mathbf{n}+(f-\bar{f})\quad\text{on}\,\ \partial\Omega,
\end{aligned}
\right.
\end{equation}
where $\bar{\mathcal{M}}=(\bar{\mathcal{M}}_{ij})$ is a  matrix
with
$\bar{\mathcal{M}}_{ij}=\dfrac{\bar{\rho}}{c^{2}(\bar{\rho})}(c^{2}(\bar{\rho})\delta_{ij}-\partial_i\bar{\phi}\partial_j\bar{\phi})$,
and
$$\mathcal{F}(x, B, W, q)= -H(B-\dfrac{1}{2}|\nabla\bar\phi+q+W|^2)(\nabla \bar\phi +q)+ H(\bar B-\dfrac{1}{2}|\nabla\bar\phi|^2)\nabla\bar\phi +  \bar{\mathcal{M}}q.$$

We prove the existence of solution of (\ref{5.4psystem}) by the contraction mapping theorem.

Given $\Psi\in \mathcal{K}=\{\Psi\in \mcC^{3,\alpha}(\bar{\Omega}): \|\nabla\Psi\|_{C^{2,\alpha}(\bar{\Omega})}\leq\theta(\leq\sigma_0/2)~\text{and}~\partial_i\Psi=\partial_{iii}\Psi=0~\text{on}~\Gamma_i\}$, it is easy to see that the problem
\begin{equation}\label{Add1}
\left\{
\begin{aligned}
&\nabla\cdot(\bar{\mathcal{M}}\nabla\Phi)=\nabla\cdot\mathcal{F}(x,B,W,\nabla\Psi)~\quad \text{in}\,\,\Omega,\\
&\bar{\mathcal{M}}\nabla\Phi\cdot\mathbf{n}=\mathcal{F}(x,B,W,\nabla\Psi)\cdot\mathbf{n}+f-\bar{f}\quad \text{on}\,\,\partial\Omega.
\end{aligned}
\right.
\end{equation}
has a solution $\Phi\in \mcA$. By the Schauder estimate, it is easy to see that $\Phi \in \mathcal{C}^{3, \alpha}(\Omega)$. Furthermore, the straightforward calculations give
\begin{equation*}
\mathcal{F}\cdot\mathbf{n}=f-\bar f=0~\text{on}~\Gamma
\end{equation*}
and
\begin{equation}
\begin{aligned}
&\partial_{ii}\mathcal{F}_i=0,\quad\partial_i\mathcal{F}_j=\mathbf{0}\quad\text{on}~\Gamma_i (i\neq j),\,\, \text{and}\,\,\frac{\partial f}{\partial \nu}=0\quad\text{on}~~\partial\Gamma_-
\end{aligned}
\end{equation}
with $\nu$ the outer normal vector of $\partial\Gamma_-$. These, together with $\partial_i\bar\phi=\partial_{iii}\bar\phi=0$ and Lemma \ref{appendix1} in Appendix \ref{appendixreflection}, imply that
$\Phi\in C^{3, \alpha}(\bar\Omega)$ and satisfies $\partial_{iii}\Phi=0$ on $\Gamma_i$ and
\begin{equation}\label{5.4nonlinearestimate}
\|\nabla\Phi\|_{{2,\alpha;\Omega}}\leq C_{1}
(\|\mathcal{F}(x,B,W,\nabla\Psi)\|_{{2,\alpha; \Omega}}+\|f-\bar{f}\|_{{2,\alpha;\Omega}})\leq C_{1}(\delta+\sigma_2+\theta^2).
\end{equation}
If $\delta_1+\sigma_2$ and $\theta$ are sufficiently small, then $\Phi \in \mathcal{K}$.

Furthermore, for $\Psi_m\in\mathcal{K}(m=1,2)$, let $\Phi_m=\mathcal{J}\Psi_m$ be the associated solution of \eqref{Add1}. Then the difference $\Phi_1-\Phi_2$ satisfies
\begin{equation}\label{add1}
\left\{
\begin{aligned}
&\Div(\bar{\mathcal{M}}\nabla(\Phi_1-\Phi_2))=\Div(\mathcal{F}(x,B,W,\nabla\Psi_1)-\mathcal{F}(x,B,W,\nabla\Psi_2))\quad&&\text{in}\,\ \Omega,\\
&\bar{\mathcal{M}}\nabla(\Phi_1-\Phi_2)\cdot\mathbf{n}=(\mathcal{F}(x,B,W,\nabla\Psi_1)-\mathcal{F}(x,B,W,\nabla\Psi_2))\cdot\mathbf{n}\quad&&\text{on}\,\ \partial\Omega.
\end{aligned}
\right.
\end{equation}
Applying Theorem \ref{appendix1} again gives
$$\|\nabla(\Phi_1-\Phi_2)\|_{{2,\alpha};\bar{\Omega}}\leq C\|\mathcal{F}(x,B,W,\nabla\Psi_1)-\mathcal{F}(x,B,W,\nabla\Psi_2)\|_{{2,\alpha};\bar{\Omega}}.$$
With the aid of the straightforward computations, one has
\begin{equation*}
\begin{array}{llllll}
\|\mathcal{F}(x,B,W,\nabla\Psi_1)-\mathcal{F}(x,B,W,\nabla\Psi_2)\|_{{2,\alpha};\bar{\Omega}}\leq C_{1}(\theta+\sigma_2)\|\nabla(\Psi_1-\Psi_2)\|_{{2,\alpha};\bar{\Omega}}.
\end{array}
\end{equation*}
Thus
$$\|\nabla(\Phi_1-\Phi_2)\|_{{2,\alpha};\bar{\Omega}}\leq C_1(\theta+\sigma_2)\|\nabla(\Psi_1-\Psi_2)\|_{{2,\alpha};\bar{\Omega}}.$$

One may choose $\theta$ and $\sigma_2$ small enough such that
\begin{equation*}
\begin{array}{llllll}
\|\mathcal{F}(x,B,W,\nabla\Psi_1)-\mathcal{F}(x,B,W,\nabla\Psi_2)\|_{{2,\alpha};\bar{\Omega}}\leq \frac{1}{2}\|\nabla(\Psi_1-\Psi_2)\|_{{2,\alpha};\bar{\Omega}}.
\end{array}
\end{equation*}
Hence (\ref{5.4nonlinear}) has a unique solution
$\phi\in \mathcal{K}$  by the contraction mapping theorem.

Moreover, it follows from (\ref{5.4nonlinearestimate}) that
\begin{equation*}
\|\nabla(\phi-\bar{\phi})\|_{{2,\alpha};\bar{\Omega}}\leq C_{1}(\|W\|_{{2,\alpha};\bar{\Omega}}+\|B-\bar{B}\|_{{2,\alpha};\bar{\Omega}}+\|f-\bar{f}\|_{{2,\alpha};\partial\Omega}+
\|\nabla(\phi-\bar{\phi})\|^2_{{2,\alpha};\bar{\Omega}}).
\end{equation*}
Therefore,
\begin{equation}
\|\nabla(\phi-\bar{\phi})\|_{{2,\alpha};\bar{\Omega}}\leq
C_1(\|W\|_{{2,\alpha};\bar{\Omega}}+\|B-\bar{B}\|_{{2,\alpha};\bar{\Omega}}+\|f-\bar{f}\|_{{2,\alpha};\partial\Omega}).
\end{equation}

On the other hand, suppose that $\phi$ and $\tilde{\phi}$ are solutions associated with $(W,B)$ and $(\tilde{W},\tilde{B})$, respectively.
Then,  $\phi-\tilde{\phi}$ satisfies the system,
\begin{equation*}
\left\{
\begin{aligned}
&\Div(\bar{\mathcal{M}}\nabla(\phi-\tilde{\phi}))=
\Div(\mathcal{F}(x,B,W,\nabla(\phi-\bar{\phi}))-\mathcal{F}(x,\tilde B,\tilde W,\nabla(\tilde{\phi}-\bar{\phi})))\quad&\text{in}\quad\Omega,\\
&\bar{\mathcal{M}}\nabla(\phi-\tilde{\phi})\cdot\mathbf{n}=(\mathcal{F}(x,B,W,\nabla(\phi-\bar{\phi}))-\mathcal{F}(x,\tilde B,\tilde W,\nabla(\tilde\phi-\bar{\phi})))\cdot\mathbf{n}\quad&\text{on}\quad\partial\Omega.
\end{aligned}
\right.
\end{equation*}
Using Theorem \ref{appendix1} again gives
\begin{equation}\label{5.4contraction}
\|\nabla(\phi-\tilde{\phi})\|_{{1,\alpha};\bar{\Omega}}\leq C(\|(B-\tilde{B}\|_{{1,\alpha};\bar{\Omega}}+\|W-\tilde{W})\|_{{1,\alpha};\bar{\Omega}}).
\end{equation}
This complete the proof of the Lemma \ref{5.4nonlinearlemma}.

\appendix
\section{Analysis of streamlines}\label{appendixstream}
In this appendix, we give the proof for Lemma \ref{lem1}.

(1) Let us rewrite (\ref{characteristicline}) into the integral form,
\begin{equation*}
\begin{aligned}
&X_i(s;x)-x_i=\int^{s}_{x_1}U_i(\tau,X_2(\tau;x),X_3(\tau;x))d\tau;\\
&\frac{\partial X_i}{\partial x_j}(s;x)-\delta_{ij}=\sum_{l=2}^{3}\int^{s}_{x_1}\dfrac{\partial U_i}{\partial x_l}(\tau,X_2(\tau;x),X_3(\tau;x))\dfrac{\partial X_l}{\partial x_j}(\tau;x)d\tau-\delta_{1j}U_i(x)
\end{aligned}
\end{equation*}
for $x=(x_1,x_2,x_3)$.
In particular, $\dfrac{\partial X_i}{\partial x_j}(x_1;x)=\delta_{ij}-\delta_{1j}U_i(x)$.
\begin{equation}
\begin{array}{lll}
\dfrac{\partial^2 X_i}{\partial x_jx_k}(s;x)&=&\sum_{l,m=2}^{3}\int^{s}_{x_1}\dfrac{\partial^2 U_i}{\partial x_l\partial x_m}(\tau,X_2(\tau;x),X_3(\tau;x))\dfrac{\partial X_l}{\partial x_j}(\tau)\dfrac{\partial X_m}
{\partial x_k}(\tau)\\
&+&\sum_{l=2}^{3}\dfrac{\partial U_i}{\partial x_l}(\tau,X_2(\tau),X_3(\tau))\dfrac{\partial^2 X_l}{\partial x_j\partial x_k}(\tau)d\tau+\delta_{1j}\dfrac{\partial U_i}{\partial x_k}(x)\\
&-&\sum_{l=2}^{3}\delta_{1k}\dfrac{\partial U_i}{\partial x_l}(x)(\delta_{lj}-\delta_{1j}U_l(x)).
\end{array}
\end{equation}
The the direction computations give the estimate \eqref{tranlabel1}.

(2) Now the difference equation can be written as follow,
\begin{equation*}
\begin{aligned}
&X_i(s)-\tilde{X}_i(s)=\int^{s}_{x_1}U_i(\tau,X_2(\tau),X_3(\tau))-\tilde{U}_i(\tau,\tilde{X}_2(\tau),\tilde{X}_3(\tau))d\tau,\\
&\frac{\partial(X_i(s)-\tilde{X}_i(s))}{\partial x_j}=\int^{s}_{x_1}\frac{\partial U_i}{\partial x_l}\frac{\partial X_l}{\partial x_j}(\tau)-\frac{\partial\tilde{U}_i}{\partial x_l}\frac{\partial\tilde{X}_l}{\partial x_j}(\tau)d\tau-\delta_{1j}(U_i(x)-\tilde{U}_i(x)).
\end{aligned}
\end{equation*}
The direct computation gives the estimate \eqref{tranlabel2}.

(3) To see this, we need to study the system for $\partial_2X_3$ and $\partial_3X_2$.
It follows from the condition $\mathbf{u}\cdot\mathbf{n}=0, (\nabla\times\mathbf{u})\times\mathbf{n}=\mathbf{0}$ on $\Gamma$ that
the function $U_i=U_i(x_1,x_2,x_3)$ satisfies $\dfrac{\partial U_3}{\partial x_2}=\dfrac{\partial U_2}{\partial x_3}=0$ on $\Gamma_2$. Therefore, for $x\in \Gamma_2$, one has
\begin{equation*}
\left\{
\begin{aligned}
&\dfrac{d(\partial_2X_3)}{ds}(s; x)=\dfrac{\partial U_3}{\partial X_3}(x_1, X_2(s; x), X_3(s;x))\partial_2X_3(s; x)=0,\\
&\partial_2X_3(x_1; x)=0.
\end{aligned}
\right.
\end{equation*}
Hence $\partial_2X_3(s; x)\equiv 0$ for $x\in \Gamma_2$. Similarly, we can prove that $\partial_3X_2(s; x)\equiv0$ for $x\in \Gamma_3$.

This finishes the proof of Lemma \ref{lem1}. {{\hfill$\Box$}

 \section{Reflection technique for linear elliptic equations}\label{appendixreflection}
In this appendix, we prove the regularity for the boundary value problem for the elliptic equation near the edges and corners of the boundary when some compatibility conditions are satisfied.

Consider the elliptic equation with cornormal boundary condition, i.e.,
\begin{equation}\label{2reflection}
\left\{
\begin{aligned}
&\Div~ (\mathcal{A}\nabla\varphi)=\Div~\mathcal{F}\quad&\text{in}\quad\Omega=[0,L]\times[0,1]^2,\\
&\mathcal{A}\nabla\varphi\cdot\mathbf{n}=\mathcal{F}\cdot\mathbf{n}+f\quad&\text{on}\quad\partial\Omega,
\end{aligned}
\right.
\end{equation}
where $\mathcal{A}=[a_{ij}]_{3\times3}$, $\mathcal{F}=(\mathcal{F}_1,\mathcal{F}_2,\mathcal{F}_3)$, and $f$ are positive definite matrix function, vector function, and scalar function, respectively.
Then we have the following global regularity results.
\begin{thm}\label{appendix1}
Suppose that
\begin{equation}\label{cond1}
\mathcal{F}\cdot\mathbf{n}=f=0~\text{on}~\Gamma
\end{equation}
and
\begin{equation}\label{2c}
\begin{aligned}
&\partial_{ii}\mathcal{F}_i=0,\quad\partial_i\mathcal{F}_j=\mathbf{0}\quad\text{on}~\Gamma_i (i\neq j),\,\, \text{and}\,\,\frac{\partial f}{\partial \nu}=0\quad\text{on}~~\partial\Gamma_-
\end{aligned}
\end{equation}
with $\nu$ the outer normal vector of $\partial\Gamma_-$.
Let
$\varphi$ be a $C^{3,\alpha}({\Omega})$ solution to (\ref{2reflection}).
\begin{enumerate}
\item[(i)] If  $a_{ij}=A(|\nabla\phi|^2)\delta_{ij}-B(|\nabla\phi|^2)\partial_i\phi\partial_j\phi$  with $\phi\in C^{3, \alpha}(\bar \Omega)$ satisfying
\begin{equation}\label{cond2}
\nabla\phi \cdot \mathbf{n}=0\quad \text{on}\,\, \Gamma\quad \text{and}\quad  \partial_{iii}\phi=0 \,\, \text{on}\,\, \Gamma_i.
\end{equation}
Then $\varphi\in C^{3, \alpha}(\bar \Omega)$ and satisfies $\partial_{iii}\varphi=0$ on $\Gamma_i$ and
\begin{equation}\label{gr}
\|\nabla\varphi\|_{C^{2,\alpha}(\bar{\Omega})}\leq C(\|\mathcal{F}\|_{C^{2,\alpha}(\bar{\Omega})}+\|f\|_{C^{2,\alpha}(\partial\Omega)}),
\end{equation}
where $C$ depends on $\|\nabla\phi\|_{C^{2,\alpha}}$ and the elliptic constants.
\item[(ii)] If $a_{ij}=\lambda\delta_{ij}$ with $\lambda$ satisfying $\nabla \lambda \cdot \mathbf{n}=0$ on $\Gamma$, then $\varphi$ also belongs to $C^{3, \alpha}(\bar\Omega)$ and satisfies \eqref{gr} with $C$ depending on $\|\lambda\|_{C^{2,\alpha}}$ and the elliptic constants..
\end{enumerate}
\end{thm}
\begin{pf}
(i) Define the extended functions for $f,\phi,\varphi$ and $\mathcal{F}=(F_1,F_2,F_3)^T$ as follows,
{\small\begin{equation*}
\begin{aligned}
&\tilde{\phi}(x_1,x_2,x_3)=\phi(x_1,|x_2-2k_2|,|x_3-2k_3|),\tilde{\varphi}(x_1,x_2,x_3)=\varphi(x_1,|x_2-2k_2|,|x_3-2k_3|),\\
&\tilde{f}(x_2,x_3)=f(|x_2-2k_2|,|x_3-2k_3|), \tilde{F}_1(x_1,x_2,x_3)=F_1(x_1,|x_2-2k_2|,|x_3-2k_3|),\\
&\tilde{F}_i(x_1,x_2,x_3)=-sgn(x_i-2k_i)F_i(x_1,|x_2-2k_2|,|x_3-2k_3|)\quad i=2,3,\\
&\tilde{\mathcal{F}}=(\tilde{F}_1,\tilde{F}_2,\tilde{F}_3)\quad \text{for}\,\ |x_i-2k_i|\leq 1, k_i=0,\pm1,\pm1\cdots,
\end{aligned}
\end{equation*}}
where $x=(x_1, x_2, x_3)$ satisfies $x_2\in (2k_2-1, 2k_2+1)$ and $x_3\in (3k_3-1, 2k_3+1)$ with $(k_2, k_3)\in \mathbb{Z}^2$. It is easy to see that
the compatibility conditions \eqref{cond1}-(\ref{cond2}) imply that $\tilde{f}\in C^{2,\alpha}(\mathbb{R}^2)$, $\nabla\tilde{\phi},\tilde{\mathcal{F}}\in C^{2,\alpha}([0,L]\times\mathbb{R}^2)$. Furthermore, $\tilde{\varphi}$ is a solution to the following boundary value problem,
\begin{equation*}
\left\{
\begin{aligned}
&\Div~(A(\nabla\tilde{\phi})\nabla\tilde{\varphi})=\Div~\tilde{\mathcal{F}}\quad&\text{in}\quad[0,L]\times\mathbb{R}^2,\\
&A(\nabla\tilde{\phi})\nabla\tilde{\varphi}\cdot\mathbf{n}=\tilde{\mathcal{F}}\cdot\mathbf{n}+\tilde{f}\quad&\text{on}\quad x_1=0,L.
\end{aligned}
\right.
\end{equation*}
Note that \eqref{cond1} implies that the problem \eqref{2reflection} has only one solution in $\mcA$.
Hence, applying the Schauder estimates for $\tilde{\varphi}$ in $[0,L]\times\mathbb{R}^2$ yields that
\begin{equation}\label{estgr}
\|\nabla\tilde{\varphi}\|_{{2,\alpha;}[0,L]\times\mathbb{R}^2}\leq C(\|\tilde{\mathcal{F}}\|_{{2,\alpha;}[0,L]\times\mathbb{R}^2}+\|\tilde{f}\|_{{2,\alpha;}\mathbb{R}^2}).
\end{equation}
This gives \eqref{gr}.

(ii) If we extend $\tilde\varphi$, $\tilde f$ and $\tilde{\mathcal{F}}$ as above and extend $\lambda$ as follows
\begin{equation*}
\begin{aligned}
&\tilde{\lambda}(x_1,x_2,x_3)=\lambda(x_1,|x_2-2k_2|,|x_3-2k_3|),
\end{aligned}
\end{equation*}
then one also has the global regularity estimate \eqref{estgr} which gives \eqref{gr}.
\end{pf}

%\newpage

\bigskip
{\bf Acknowledgement.}    The work is part of the first author's Ph.D thesis at The Chinese University of Hong Kong under the supervision of Professor Zhouping Xin. Part of the work was done when the second author was visiting The Institute of Mathematical Sciences, The Chinese University of Hong Kong. Both of the authors thank Professor Zhouping Xin for helpful discussions.
The research was support in part by NSFC grants 11241001 and 11201297, Shanghai Chenguang program and Shanghai Pujiang  program 12PJ1405200.

\end{document}